\documentclass[10pt]{article}
\usepackage{amsmath,amssymb}
\usepackage{booktabs}
\usepackage[margin=2.5cm]{geometry}
\usepackage[hidelinks]{hyperref}

\title{The mean orders of even and odd permutations}
\author{Manoj Kumar Singh\\[2pt]
\small Department of Mathematics and Statistics, St.~Xavier's College (Autonomous), Ranchi 834001, India\\
\small \texttt{manojkumarsingh@sxcran.org}}
\date{}

\begin{document}
\maketitle

\begin{abstract}
We compare the sums of the orders of the even and odd permutations of a given degree. For all
large degrees the odd permutations have the larger sum. The excess equals the sum of orders of the
permutations of odd order up to a relative error decaying faster than any power of the degree.
Following Goh and Schmutz, we determine the logarithm of this sum over the group order
asymptotically, with an explicit constant. Hence the ratio of the sums of element orders of the
symmetric and alternating groups tends to two from above. The proof rests on an exact identity for
generating functions. It isolates cycles of length a large prime times a small factor. Inclusion
and exclusion over small odd primes enlarges this family. Each Euler factor is bounded by its
maximum modulus. All constants are certified in ball arithmetic. Exact values up to degree 150
suggest that the odd permutations have the larger sum from degree four on.
\end{abstract}

\noindent\textit{Mathematics Subject Classification (2020).} 60C05, 05A16, 20B30, 11N37.

\noindent\textit{Keywords.} Order of a random permutation; sign character; symmetric group;
alternating group; sum of element orders.

\section{Introduction}

For a permutation $\sigma$ let $o(\sigma)$ denote its order. Put
\[
T_n=\sum_{\sigma\in S_n}o(\sigma),
\]
so that $T_n/n!$ is the mean order of a random permutation of $n$ points. Put
$s_n=\sqrt{n/\log n}$. Erd\H{o}s and Tur\'an \cite{ET1967} proved that $\log o(\sigma)$ is
asymptotically normal with mean $\tfrac12\log^2n$ and variance $\tfrac13\log^3n$. The mean order
is much larger than the typical order. It is carried by a thin set of permutations. Erd\H{o}s
and Tur\'an \cite{ET1968} raised the problem of estimating the mean order. They determined that
$\log(T_n/n!)=O(s_n)$, as reported in \cite{GohSchmutz1991}. A proof appears in
\cite{Schmutz1989}. Goh and Schmutz
\cite{GohSchmutz1991} proved that $\log(T_n/n!)\sim c\,s_n$ with an explicit constant $c$. They
remark that an asymptotic formula for $T_n/n!$ itself is probably out of reach. Stong
\cite{Stong1998} proved
\begin{equation}\label{eq:GS}
\log\frac{T_n}{n!}=c\,s_n+O\Bigl(\frac{\sqrt n\,\log\log n}{\log n}\Bigr),\qquad c=2\sqrt{2A},\qquad
A=\int_0^\infty\log\bigl(1-\log(1-e^{-r})\bigr)\,dr .
\end{equation}
This is stated in the introduction of \cite{Stong1998}. There $A$ is written as the integral of
$re^{-r}/((1-e^{-r})(1-\log(1-e^{-r})))$. Integration by parts gives the form above, because
$r\log(1-\log(1-e^{-r}))$ tends to $0$ at both ends. Numerically $A=1.1178\ldots$ and
$c=2.9904\ldots$ The error term in \eqref{eq:GS} is $o(s_n)$.

The order of a random permutation has been studied in many other settings. Beker \cite{Beker2025}
found the most probable order of a random permutation of $n$ points for large $n$. In
\cite{Beker2026} he determined the asymptotic behaviour of the R\'enyi entropies of the order.
Acan, Burnette, Eberhard, Schmutz and Thomas \cite{ABEST2021} studied the probability that two random
permutations have the same order. Ford \cite{Ford2022} gave a toolkit of estimates for the cycle
type of a random permutation. Storm and Zeindler \cite{StormZeindler} gave large deviation estimates
for the order of random permutations with polynomial cycle weights. Yakymiv proved limit theorems
for the logarithm of the order of random $\mathcal A$-permutations \cite{Yakymiv2010}, whose cycle
lengths lie in a fixed set $\mathcal A$. He also treated permutations with cycle weights \cite{Yakymiv2018}. These works
concern the typical order or its logarithm. The mean order depends on a thin set of permutations.
We know of no estimate of the mean order in these settings. The mean order also enters
applications. Bavisetty et al. \cite{Bavisetty2026} used the theorem of Goh and Schmutz in a study
of the attractors of critical Boolean networks. The permutations of odd order are the
$\mathcal A$-permutations with $\mathcal A$ the set of odd numbers. Their exponential generating
function is $\sqrt{(1+x)/(1-x)}$. So their proportion in $S_n$ is asymptotic to $\sqrt{2/(\pi n)}$.
Hence Theorem 1.1(a) below also gives the logarithm of their mean order.

This paper asks how the mean order is shared between the even and the odd permutations. Put
\[
D_n=\sum_{\sigma\in S_n}\operatorname{sgn}(\sigma)\,o(\sigma).
\]
The sum of the orders over the alternating group $A_n$ is $(T_n+D_n)/2$. The sum over the odd
permutations is $(T_n-D_n)/2$. The signed sum $D_n$ corresponds to the cycle weights
$(-1)^{\ell-1}$. Let $c(\sigma)$ be the number of cycles of $\sigma$. Then
$\operatorname{sgn}(\sigma)=(-1)^{n-c(\sigma)}$. So $D_n=(-1)^n\sum_\sigma(-1)^{c(\sigma)}o(\sigma)$. This is
the Ewens-type weighting $\theta^{c(\sigma)}$ with $\theta=-1$. These weights do not define a
probability measure. We found no earlier
estimate of $D_n$.

The answer involves the permutations of odd order. Put
\[
V_n=\sum_{\substack{\sigma\in S_n\\ o(\sigma)\ \mathrm{odd}}}o(\sigma).
\]

\medskip\noindent\textbf{Theorem 1.1.} \textit{{\rm(a)} As $n\to\infty$,}
\[
\log\frac{V_n}{n!}\sim c_2\,s_n,\qquad c_2=2\sqrt{2A_2},\qquad
A_2=\int_0^\infty\log\Bigl(1+\frac12\log\frac{1+e^{-r}}{1-e^{-r}}\Bigr)dr .
\]
\textit{{\rm(b)} There is a constant $\delta>0$ such that, as $n\to\infty$,}
\[
D_n=-V_n\bigl(1+O(e^{-\delta s_n})\bigr).
\]
\textit{One can take $\delta=0.07$.}

\medskip
Numerically $A_2=0.9356\ldots$ and $c_2=2.7359\ldots$ So the excess of the odd permutations over
the even ones is asymptotically the sum of the orders of the permutations of odd order.

\medskip\noindent\textbf{Corollary 1.2.} \textit{{\rm(a)} We have $D_n<0$ for all sufficiently
large $n$. So the odd permutations of $\{1,\dots,n\}$ have the larger sum of orders.}

\textit{{\rm(b)} We have $\log(|D_n|/T_n)\sim-(c-c_2)s_n$, where $c-c_2=0.2545\ldots$}

\medskip
The question arose from the function $\psi(G)=\sum_{g\in G}o(g)$ of a finite group $G$. It was
introduced by Amiri, Jafarian Amiri and Isaacs \cite{AAI2009}. Recent work on $\psi$ includes
\cite{Amiri2025,CDP2026,IorioTrombetti2026,Tarnauceanu2026}. The quotient $\psi(G)/|G|$ is the average order of $G$.
Herzog, Longobardi and Maj \cite{HLM2022} proved that $\psi(G)/|G|<\psi(A_5)/60$ forces $G$ to be
solvable. In this notation $T_n=\psi(S_n)$. So $T_n/n!$ is the average order of $S_n$. Goh and
Schmutz determined its logarithm asymptotically. The works above do not separate the two cosets of
$A_n$. We have $(T_n+D_n)/2=\psi(A_n)$. Hence
\begin{equation}\label{eq:ratio}
\frac{\psi(S_n)}{\psi(A_n)}=\frac{2}{1+D_n/T_n}.
\end{equation}
The inequality $3\psi(A_n)>\psi(S_n)$ was proposed for $n\ge3$ in \cite{Singh2026}. That preprint
reduces it to a statement on the number of odd square roots of the even permutations. It verifies
the inequality for $n\le60$. It also asks whether the ratio $\psi(S_n)/\psi(A_n)$ is largest at
$n=13$. The inequality says that $D_n>-T_n/3$. At $n=2$ equality holds.

\medskip\noindent\textbf{Corollary 1.3.} \textit{We have
$\psi(S_n)/\psi(A_n)=2+(2+o(1))V_n/T_n$. So $\psi(S_n)/\psi(A_n)>2$ for all sufficiently large
$n$.}

\medskip
The limit of this ratio needs less than Theorem 1.1(b). Proposition 6.3 states that
$|D_n|\le T_ne^{-0.19s_n}$ for all sufficiently large $n$. Its proof applies the bounds of Sections 4
and 5 directly. It does not use the inclusion and exclusion of Lemma 6.1. Hence $D_n=o(T_n)$. So the
even and the odd permutations have asymptotically equal mean order. Also $\psi(S_n)/\psi(A_n)\to2$
and $3\psi(A_n)>\psi(S_n)$ for all sufficiently large $n$. Theorem 1.1(b) is needed for three finer
statements. These are the sign of $D_n$, the exact exponent in Corollary 1.2(b) and the approach of
the ratio to $2$ from above.

For $n\le150$ we computed $T_n$, $D_n$ and $V_n$ exactly.

\medskip\noindent\textbf{Proposition 1.4.} \textit{{\rm(a)} For $3\le n\le150$ we have
$3\psi(A_n)>\psi(S_n)$. Among these $n$ the ratio $\psi(S_n)/\psi(A_n)$ is largest at $n=13$,
where it equals $2.56117\ldots$}

\textit{{\rm(b)} For $n=2$ and for $4\le n\le150$ we have $D_n<0$.}

\medskip
Theorem 1.1, Corollary 1.2, Corollary 1.3 and Proposition 6.3 are asymptotic statements. Section 7
discusses their thresholds. Conjecture 7.1 states what we expect for every $n$.

We now describe the proof of Theorem 1.1. The method of Goh and Schmutz bounds $T_n/n!$ through
the product $\prod_p(1-\log(1-e^{-pt}))$. Here $-\log(1-e^{-pt})=\sum_jj^{-1}e^{-jpt}$ collects the
cycles of length $jp$. In the signed sum a cycle of length $jp$ with $p$ odd carries the sign
$(-1)^{j-1}$. If these signs were free, the factor would become $1+\log(1+e^{-pt})$. A heuristic
count with free signs then gives the exponent $2\sqrt{2A_-}=2.3775\ldots$ with
\[
A_-=\int_0^\infty\log\bigl(1+\log(1+e^{-r})\bigr)\,dr=0.7066\ldots .
\]
This is smaller than $c_2$ because $\log(1+u)<\tfrac12\log\frac{1+u}{1-u}$ for $0<u<1$. This is the
heuristic reason behind Theorem 1.1(b). The signs are not free. An even multiplier $j$ also
doubles the order. Lemma 3.1 measures this effect exactly. It states that $D_n=D_n^*-V_n$, where
$D_n^*$ is the sum of $\operatorname{sgn}(\sigma)\operatorname{lcm}(2,o(\sigma))$. The sum $V_n$ has
positive terms. Its size follows from the method of Goh and Schmutz for partitions into odd parts
(Section 3). Part (a) of Theorem 1.1 is thus an adaptation of their method. The new content of the
paper is part (b).

The signed sum $D_n^*$ is smaller. To prove this we write it as a combination of further signed
sums by inclusion and exclusion over a set $Q$ of small odd primes (Lemma 6.1). The set $Q=\{3\}$
gives two sums. They already suffice. We use $Q=\{3,5\}$, which gives four sums and a better rate.
For each of them an exact identity separates the cycles whose lengths
are a large prime times a small divisor (Section 4). The other cycles give a positive series. We
bound it by the decomposition of Goh and Schmutz. The separated cycles give a product of entire
functions, one for each large prime. On a circle we bound each factor by its maximum modulus
(Section 5). Cauchy's inequality then bounds the signed sum. The resulting exponents have the form
$2\sqrt{2(I+\nu)}$ with explicit integrals $I$ and $\nu$. The largest value of $I+\nu$ is
numerically $0.8834\ldots$ It exceeds $A_-$ because of the order-doubling effect and the bounds used
for the non-free cycles. Certified computations show that it is smaller than $A_2$ (Section 6). The
largest value is at most $0.8839$, while $A_2\ge0.9356$.
Section 7 reports the exact computations.

The proof of part (b) adds three ingredients to the method of Goh and Schmutz. The first is the
reduction $D_n=D_n^*-V_n$ of Lemma 3.1. It isolates the permutations of odd order as the main term.
The second is the identity of Proposition 4.1 together with the maximum-modulus bound of Section 5.
They replace the positivity on which the method of Goh and Schmutz rests. The third is the inclusion
and exclusion of Lemma 6.1. It makes the small multipliers free. Section 2 and the rest of
Section 3 follow \cite{GohSchmutz1991} closely. The new arguments begin in Section 4. The bounds of Sections 4 and 5
hold for every fixed pair $(m,e)$. So they apply to other signed sums of this type.

\section{Preliminaries}

\subsection{Notation}

Throughout, $p$ denotes a prime and $\pi(u)$ the number of primes up to $u$. For a cycle length
$\ell$ put $\varepsilon_\ell=(-1)^{\ell-1}$. A permutation with $a_\ell$ cycles of length $\ell$ has
sign $\prod_\ell\varepsilon_\ell^{a_\ell}$. Its order is the least common multiple of its cycle
lengths. When $t\to0^+$ we put $L=\log(1/t)$.

\subsection{The decomposition of Goh and Schmutz}

The following construction goes back to Schmutz \cite{Schmutz1989}. We use it in the form given by
Goh and Schmutz \cite[Section 3]{GohSchmutz1991}. Let
$\lambda$ be a finite multiset of positive integers. Let
$\operatorname{lcm}(\lambda)=p_1^{e_1}\cdots p_s^{e_s}$ with $p_1<\dots<p_s$. Choose elements
$b_1,b_2,\dots$ of $\lambda$ as follows. At each stage take the least $k$ such that $p_k^{e_k}$
divides none of the elements already chosen. Choose the smallest element of $\lambda$ divisible by
$p_k^{e_k}$. Stop when every $p_k^{e_k}$ divides a chosen element. Assign each $p_k^{e_k}$ to the
first chosen element that it divides. Let $d_i$ be the product of the prime powers assigned to
$b_i$. Put $d_i'=b_i/d_i$. Let $\omega$ be the multiset of the remaining elements.

\medskip\noindent\textbf{Lemma 2.1.} \textit{The $d_i$ are pairwise
coprime. Each $d_i$ exceeds $1$ and divides $b_i$. Their product is $\operatorname{lcm}(\lambda)$.
The multiset $\lambda$ is determined by the set of pairs $(d_i,d_i')$ together with $\omega$.
Finally}
\[
\frac{\operatorname{lcm}(\lambda)}{\prod_{x\in\lambda}x}=\frac{1}{\prod_id_i'\prod_{y\in\omega}y}.
\]

\medskip
These properties are stated in \cite[Section 3]{GohSchmutz1991}. We have $d_i>1$ because $b_i$
was chosen for a prime power that is then assigned to it.

\subsection{Three auxiliary results}

The first result is a version of \cite[Lemma 2]{GohSchmutz1991} for a general monotone function.
We use the same splitting points. Their proof uses the prime number theorem with an error term. The
proof below needs only the asymptotic formula $\pi(r)\sim r/\log r$ and a Chebyshev bound. We include
it because the lemma is applied to four different functions.

\medskip\noindent\textbf{Lemma 2.2.} \textit{Let $g:(0,\infty)\to[0,\infty)$ be continuous and
non-increasing. Suppose that for some $B\ge1$ we have $g(r)\le B\log(2+\log(1/r))$ for
$0<r\le1$ and $g(r)\le Be^{-r}$ for $r\ge1$. Put $I_g=\int_0^\infty g$. Then, as $t\to0^+$,}
\[
\sum_pg(pt)=\frac{I_g+o(1)}{t\log(1/t)} .
\]
\textit{For every fixed $y\ge1$ the same holds for the sum over the primes $p>y$.}

\medskip\noindent\textit{Proof.} Put $\xi=1/(tL^3)$ and $\xi'=L^3/t$. Let $\mu$ be the
Lebesgue--Stieltjes measure with $\mu((a,b])=g(at)-g(bt)$. It is non-negative. Since $g$ is
continuous and tends to $0$, we have $g(pt)=\mu((p,\infty))$. By Tonelli's theorem
\begin{equation}\label{eq:stieltjes}
\sum_pg(pt)=\int_{(2,\infty)}\#\{p<r\}\,d\mu(r).
\end{equation}

We first prove the upper bound. We have $\#\{p<r\}\le\pi(r)$. On $(2,\xi]$ we use
$\pi(r)\le\pi(\xi)\le2\xi/\log\xi$ \cite{RosserSchoenfeld} and $\mu((2,\xi])\le g(2t)$. For small
$t$ we have $\log\xi\ge L/2$ and $g(2t)\le B\log(2+L)$. So this part of \eqref{eq:stieltjes} is at
most $4B\log(2+L)/(tL^4)$. On $(\xi,\infty)$ fix $\epsilon>0$. By the prime number theorem
$\pi(r)\le(1+\epsilon)r/\log r$ for $r\ge r_\epsilon$. For small $t$ we have $\xi\ge r_\epsilon$. There
$\log r\ge L-3\log L$. So this part is at most
$(1+\epsilon)(L-3\log L)^{-1}\int_{(\xi,\infty)}r\,d\mu(r)$. Writing $r=\xi+\int_\xi^rdu$ and using
Tonelli's theorem again,
\[
\int_{(\xi,\infty)}r\,d\mu(r)=\xi g(\xi t)+\int_\xi^\infty g(ut)\,du\le\xi g(\xi t)+\frac{I_g}{t}.
\]
Here $\xi t=L^{-3}\le1$. So $\xi g(\xi t)\le B\log(2+3\log L)/(tL^3)=o(1/t)$. The two parts give
$tL\sum_pg(pt)\le(1+\epsilon)(I_g+o(1))$. Now let $\epsilon\to0$.

We turn to the lower bound. For $r>2$ we have $\#\{p<r\}\ge\pi(r)-1$. By the prime number theorem
$\pi(r)\ge(1-\epsilon)r/\log r$ for $r\ge r_\epsilon$. On $(\xi,\xi']$ we have
$\log r\le L+3\log L$. So \eqref{eq:stieltjes} gives
\[
\sum_pg(pt)\ge\frac{1-\epsilon}{L+3\log L}\int_{(\xi,\xi']}r\,d\mu(r)-\mu((\xi,\xi']).
\]
Here $\mu((\xi,\xi'])\le g(\xi t)\le B\log(2+3\log L)$. This is $o(1/(tL))$. As above,
\[
\int_{(\xi,\xi']}r\,d\mu(r)=\xi\bigl(g(\xi t)-g(\xi't)\bigr)+\int_\xi^{\xi'}\bigl(g(ut)-g(\xi't)\bigr)du
\ge\frac1t\int_{L^{-3}}^{L^3}g(v)\,dv-\xi'g(\xi't).
\]
We have $\xi'g(\xi't)\le BL^3e^{-L^3}/t=o(1/t)$. Since $g$ is integrable, the integral of $g$ over
$[L^{-3},L^3]$ tends to $I_g$. So $tL\sum_pg(pt)\ge(1-\epsilon)(I_g-o(1))$. Now let
$\epsilon\to0$.

Finally the primes $p\le y$ contribute at most $\pi(y)g(2t)\le\pi(y)B\log(2+L)$. This is
$o(1/(tL))$. \hfill$\square$

\medskip\noindent\textbf{Lemma 2.3.} \textit{For $k\ge0$ put $\Lambda_k=\sum1/(y_1y_2\cdots)$, the sum
running over the partitions $y_1+y_2+\dots=k$. Then $\Lambda_k\le e^{3/2}(k+1)^2$.}

\medskip\noindent\textit{Proof.} We have $\Lambda_0=1$. Let $k\ge1$. The generating function is
$\sum_k\Lambda_kx^k=\prod_{j\ge1}(1-x^j/j)^{-1}$. Let $0<x<1$. For $j\ge2$,
\[
-\log\Bigl(1-\frac{x^j}{j}\Bigr)\le\frac{x^j}{j}+\sum_{r\ge2}\frac{1}{rj^r}
\le\frac{x^j}{j}+\frac{1}{2j(j-1)} .
\]
Summing over $j$ gives $\log\prod_j(1-x^j/j)^{-1}\le-2\log(1-x)+\tfrac12$. The coefficients are
non-negative. So $\Lambda_k\le x^{-k}e^{1/2}(1-x)^{-2}$. Take $x=k/(k+1)$. Then
$x^{-k}=(1+1/k)^k\le e$. \hfill$\square$

\medskip
Goh and Schmutz \cite{GohSchmutz1991} quote a bound $O(n)$, uniform in $k\le n$, from a theorem of Lehmer
\cite{Lehmer1972}. A polynomial bound is enough for us. The third result is a Tauberian theorem
of Hardy and Ramanujan \cite[Section 4, Theorem A]{HardyRamanujan1917}. We need the case
$\lambda_n=n$, $\alpha=1$, $\beta=1$ of their theorem. It is stated in this form in
\cite[Lemma 1]{GohSchmutz1991}.

\medskip\noindent\textbf{Lemma 2.4.} \textit{Let $a_k\ge0$ and $f(t)=\sum_ka_ke^{-kt}$. Suppose that
$\log f(t)\sim A_0/(t\log(1/t))$ as $t\to0^+$ with a constant $A_0>0$. Then
$\log\sum_{k=1}^na_k\sim2\sqrt{2A_0}\,s_n$ as $n\to\infty$.}

\section{Permutations of odd order}

Put
\[
D^{*}_n=\sum_{\sigma\in S_n}\operatorname{sgn}(\sigma)\operatorname{lcm}\bigl(2,o(\sigma)\bigr).
\]

\noindent\textbf{Lemma 3.1.} \textit{We have $D_n=D^{*}_n-V_n$.}

\medskip\noindent\textit{Proof.} If $o(\sigma)$ is even then
$\operatorname{lcm}(2,o(\sigma))=o(\sigma)$. Let $o(\sigma)$ be odd. Then every cycle of $\sigma$ has
odd length. So $\sigma$ is even. Also $\operatorname{lcm}(2,o(\sigma))-o(\sigma)=o(\sigma)$. Summing
over $\sigma$ gives $D^{*}_n-D_n=V_n$. \hfill$\square$

\medskip
Put $h(r)=\log(1-\log(1-e^{-r}))$ and
\[
h_2(r)=\log\Bigl(1+\sum_{j\ \mathrm{odd}}\frac{e^{-jr}}{j}\Bigr)
=\log\Bigl(1+\frac12\log\frac{1+e^{-r}}{1-e^{-r}}\Bigr).
\]
Then $A=\int_0^\infty h$ and $A_2=\int_0^\infty h_2$. Since $0<h_2<h$ on $(0,\infty)$, we have
$A_2<A$ and $c_2<c$. For $0<r\le1$ we have $1-e^{-r}\ge re^{-r}$. So $h(r)\le\log(2+\log(1/r))$. For
$r\ge1$ we have $h(r)\le e^{-r}/(1-e^{-r})\le2e^{-r}$. So every continuous non-increasing function
between $0$ and $h$ satisfies the hypotheses of Lemma 2.2 with $B=2$.

A permutation has odd order exactly when all its cycle lengths are odd. For a partition $\rho$
with $c_i$ parts equal to $i$ put
$w(\rho)=\operatorname{lcm}(\rho)/\prod_ic_i!\,i^{c_i}$. By Cauchy's formula for the size of a
conjugacy class,
\begin{equation}\label{eq:Vcauchy}
\frac{V_n}{n!}=\sum_\rho w(\rho),
\end{equation}
the sum running over the partitions $\rho$ of $n$ into odd parts.

\medskip\noindent\textbf{Proposition 3.2.} \textit{We have $\log(V_n/n!)\le(c_2+o(1))s_n$.}

\medskip\noindent\textit{Proof.} We adapt the upper bound of \cite[Section 3]{GohSchmutz1991} to
partitions into odd parts. We have
\[
w(\rho)\le\frac{\operatorname{lcm}(\rho)}{\prod_{x\in\rho}x} .
\]
Apply Lemma 2.1 to each $\rho$ in \eqref{eq:Vcauchy}. It gives a set $\mathcal D$ of pairs
$(d_i,d_i')$ and a multiset $\omega$. Every $d_id_i'$ is a part of $\rho$. So it is odd. For a fixed
$\mathcal D$ the multiset $\omega$ is a partition of the fixed integer $k=n-\sum dd'$. Its parts are
odd. We drop this condition, which only enlarges the sum. By Lemma 2.3 the sum of
$1/\prod_{y\in\omega}y$ over all partitions $\omega$ of $k$ is at most $e^{3/2}(n+1)^2$. Hence
\[
\frac{V_n}{n!}\le e^{3/2}(n+1)^2\sum_{\mathcal D}\prod_{(d,d')\in\mathcal D}\frac1{d'} .
\]
Here $\mathcal D$ runs over the finite sets of pairs $(d,d')$ of odd positive integers with the
$d$ pairwise coprime and greater than $1$ and with $\sum dd'\le n$.

Put $z=\sqrt n/\log^2n$. Split each $\mathcal D$ into three groups. In $\mathcal G_1$ the entry $d$
has a prime factor below $z$. In $\mathcal G_2$ the entry $d$ is a prime $p\ge z$. In
$\mathcal G_3$ the entry $d$ has at least two prime factors counted with multiplicity, all at least
$z$. Drop the
coprimality between the groups. Bound the sum of $dd'$ over each group separately by $n$. Then
the sum over $\mathcal D$ is at most a product $\Sigma_1\Sigma_2\Sigma_3$. Here $\Sigma_i$ is the
sum over all possible groups $\mathcal G_i$.

The entries $d$ are pairwise coprime. So each prime below $z$ divides at most one $d$ in
$\mathcal G_1$. Hence $\mathcal G_1$ has at most $\lfloor z\rfloor$ elements, because the map sending
$d$ to its least prime factor is injective on $\mathcal G_1$. There are at most
$(n+1)^{\lfloor z\rfloor}$ choices of its entries $d$. Each $d'$ is at most $n$. Since
$\sum_{d'\le n}1/d'\le1+\log n$, the entries $d'$ contribute at most
$(1+\log n)^{\lfloor z\rfloor}$. Hence
$\Sigma_1\le((n+1)(1+\log n))^z=e^{O(\sqrt n/\log n)}$. Every $d$ in $\mathcal G_3$ is at least
$z^2=n/\log^4n$. So $\mathcal G_3$ has at most $\log^4n$ elements. Hence
$\Sigma_3\le((n+1)(1+\log n))^{\log^4n}=e^{O(\log^5n)}$.

In $\mathcal G_2$ the primes are distinct and the entries $d'$ are odd. Thus
$\Sigma_2\le\sum_{k\le n}a_k$, where
$\sum_ka_kx^k=\prod_{p\ge z}\bigl(1+\sum_{j\ \mathrm{odd}}x^{pj}/j\bigr)$. The $a_k$ are
non-negative. So for every $t>0$,
\[
\Sigma_2\le\sum_ka_ke^{(n-k)t}\le\exp\Bigl(nt+\sum_ph_2(pt)\Bigr).
\]
By Lemma 2.2, $\sum_ph_2(pt)=(A_2+o(1))/(tL)$. Take $t=\sqrt{2A_2/(n\log n)}$. Then
$L=(\tfrac12+o(1))\log n$. So $nt=\sqrt{2A_2}\,s_n$ and
$(A_2+o(1))/(tL)=(\sqrt{2A_2}+o(1))s_n$. The other factors are $e^{o(s_n)}$, because
$\sqrt n/\log n=s_n/\sqrt{\log n}$. \hfill$\square$

\medskip\noindent\textbf{Proposition 3.3.} \textit{We have $\log(V_n/n!)\ge(c_2-o(1))s_n$.}

\medskip\noindent\textit{Proof.} We follow the lower bound of \cite[Section 2]{GohSchmutz1991}
closely. Only the completion to a partition into odd parts is new. Fix
an integer $y\ge2$. For $k\ge0$ let $\mathcal P_k$ be the set of partitions of $k$ into parts of
the form $dp$ with $d\le y$ odd and $p>y$ prime. In addition each prime larger than $y$ must
divide at most one part.
The parts of such a partition are odd, because $p>2$. They are pairwise distinct. Indeed
$d_1p_1=d_2p_2$ with $p_1\ne p_2$ would force $p_1$ to divide $d_2$, but $d_2\le y<p_1$.

Let $k\le n$ and $\lambda\in\mathcal P_k$. We extend $\lambda$ to a partition $\rho(\lambda)$ of $n$
into odd parts. If $k=n$ we add nothing. If $n-k$ is odd we add the part $n-k$. If $n-k$ is even
and positive we add the parts $n-k-1$ and $1$. The pair $(k,\rho(\lambda))$ determines $\lambda$.
So each $\rho$ arises from at most $n+1$ partitions $\lambda$. By \eqref{eq:Vcauchy},
\[
\frac{V_n}{n!}\ge\frac1{n+1}\sum_{k\le n}\sum_{\lambda\in\mathcal P_k}w(\rho(\lambda)).
\]
Let $\lambda$ have the parts $d_1p_1,\dots,d_qp_q$. The primes $p_i$ are distinct. So
$\operatorname{lcm}(\rho(\lambda))\ge p_1\cdots p_q$. The added parts have product at most $n$. The
parts of $\lambda$ are distinct and exceed $1$. If $n-k=2$ the added parts are $1$ and $1$. Then the
part $1$ occurs twice and every other part occurs once. In all other cases the added parts are
distinct. The part $1$ is not a part of $\lambda$. So at most one added part equals a part of
$\lambda$. In both cases at most one part size of $\rho(\lambda)$ occurs twice and none occurs three
times. Hence $\prod_ic_i!\le2$ for $\rho(\lambda)$. This gives
$w(\rho(\lambda))\ge1/(2n\,d_1\cdots d_q)$. Therefore
\[
\frac{V_n}{n!}\ge\frac1{2n(n+1)}\sum_{k\le n}b_k,\qquad
\sum_kb_kx^k=\prod_{p>y}\Bigl(1+\sum_{\substack{d\le y\\ d\ \mathrm{odd}}}\frac{x^{dp}}d\Bigr).
\]

Put $h_{2,y}(r)=\log\bigl(1+\sum_{d\le y,\,d\ \mathrm{odd}}e^{-dr}/d\bigr)$ and
$f(t)=\sum_kb_ke^{-kt}$. Then $\log f(t)=\sum_{p>y}h_{2,y}(pt)$. The function $h_{2,y}$ is
continuous and non-increasing. It lies between $0$ and $h$. By Lemma 2.2,
$\log f(t)\sim A_{2,y}/(tL)$ with $A_{2,y}=\int_0^\infty h_{2,y}>0$. By Lemma 2.4,
$\log\sum_{k=1}^nb_k\sim2\sqrt{2A_{2,y}}\,s_n$. So
\[
\liminf_{n\to\infty}\frac{\log(V_n/n!)}{s_n}\ge2\sqrt{2A_{2,y}} .
\]
As $y\to\infty$ the functions $h_{2,y}$ increase to $h_2$. So $A_{2,y}\to A_2$ by monotone
convergence. \hfill$\square$

\medskip
Propositions 3.2 and 3.3 give Theorem 1.1(a).

\section{Signed sums: free cycles and the non-free part}

Fix integers $m\ge1$ and $e\ge1$. Put
\[
S_n(m,e)=\sum_{\substack{\sigma\in S_n\\ \gcd(o(\sigma),e)=1}}\operatorname{sgn}(\sigma)
\operatorname{lcm}\bigl(m,o(\sigma)\bigr),\qquad S_0(m,e)=m .
\]
So $D_n=S_n(1,1)$ and $D_n^*=S_n(2,1)$. The condition $\gcd(o(\sigma),e)=1$ holds exactly when
every cycle length of $\sigma$ is coprime to $e$. Call a length coprime to $e$ \emph{admissible}.
Put
\[
J=\{j\ge1:\ j\mid m,\ \gcd(j,e)=1\},\qquad R=\{j\ge1:\ j\nmid m,\ \gcd(j,e)=1\},
\]
\[
G(u)=\sum_{j\in J}\frac{\varepsilon_ju^j}{j},\qquad H=\sum_{j\in J}\frac1j .
\]
Let $z>2$ exceed every prime factor of $me$. Call a prime $p\ge z$ \emph{large}. Call a length
$\ell$ \emph{free} if $\ell\in J$ or if $\ell=jp$ with $j\in J$ and $p$ large. Free lengths are
admissible. Call all other admissible lengths \emph{non-free}. For a set $X$ of lengths let
$\operatorname{pr}(X)$ be the set of large primes dividing some element of $X$. Let
$\operatorname{lcm}(m,X)$ be the least common multiple of $m$ and the elements of $X$. Put
\[
\Gamma_p(w)=1+p\bigl(e^{G(w^p)/p}-1\bigr).
\]

\medskip\noindent\textbf{Proposition 4.1.} \textit{Let $n\ge1$. Modulo $w^{n+1}$,}
\[
\sum_{k=0}^n\frac{S_k(m,e)}{k!}w^k\equiv\sum_X\operatorname{lcm}(m,X)
\prod_{\ell\in X}\bigl(e^{\varepsilon_\ell w^\ell/\ell}-1\bigr)\,e^{G(w)}
\prod_{p\in\operatorname{pr}(X)}e^{G(w^p)/p}\prod_{\substack{z\le p\le n\\ p\notin\operatorname{pr}(X)}}\Gamma_p(w),
\]
\textit{where $X$ runs over the sets of non-free lengths with sum at most $n$.}

\medskip\noindent\textit{Proof.} Let $\sigma$ be a permutation of $k$ points whose cycle lengths are
admissible. Let $X$ be the set of its non-free cycle lengths. Let $Q$ be the set of the large
primes $p\notin\operatorname{pr}(X)$ for which $\sigma$ has a cycle of length $jp$ with $j\in J$. We
claim that
\[
\operatorname{lcm}\bigl(m,o(\sigma)\bigr)=\operatorname{lcm}(m,X)\prod_{p\in Q}p .
\]
Indeed a cycle length $j\in J$ divides $m$. A cycle length $jp$ with $j\in J$ and $p$ large adds
at most the prime $p$, since $j$ divides $m$. The prime $p$ does not divide $m$. So it divides
$\operatorname{lcm}(m,X)$ exactly when $p\in\operatorname{pr}(X)$. Distinct large primes are
coprime. This proves the claim.

We now use the exponential formula. The permutations all of whose cycles have length $\ell$ have
the signed exponential generating function $e^{\varepsilon_\ell w^\ell/\ell}$. Those with at
least one cycle give $e^{\varepsilon_\ell w^\ell/\ell}-1$. The cycles with lengths in $J$ give
$e^{G(w)}$. Let $p$ be large. Since $p$ is odd, we have $\varepsilon_{jp}=\varepsilon_j$. A free length
has only one representation $jp$, because $j$ has no large prime factor. Hence the cycles with
lengths $jp$, $j\in J$, give $e^{G(w^p)/p}$. If $p\in\operatorname{pr}(X)$ these cycles do not
change the weight. If $p\notin\operatorname{pr}(X)$ the non-empty families of such cycles get the
weight $p$, which gives the factor $e^{G(w^p)/p}+(p-1)\bigl(e^{G(w^p)/p}-1\bigr)=\Gamma_p(w)$. The product of all these factors, summed over
$X$, is the generating function of $\operatorname{sgn}(\sigma)\operatorname{lcm}(m,o(\sigma))$ over
the permutations with admissible cycle lengths. The empty permutation gives the term $m$. A set
$X$ with sum above $n$ contributes only to powers above $w^n$. For $p>n$ we have
$\Gamma_p(w)\equiv1$ modulo $w^{n+1}$. \hfill$\square$

\medskip
Let $\Psi_n$ denote the right side of Proposition 4.1. It is a finite sum of finite products of
entire functions. So it is entire. By Proposition 4.1 the coefficient of $w^n$ in $\Psi_n$ is
$S_n(m,e)/n!$. For $0<x<1$ put $M_p(x)=\max_{|w|=x}|\Gamma_p(w)|$.

\medskip\noindent\textbf{Lemma 4.2.} \textit{For $0<x<1$,}
\[
\frac{|S_n(m,e)|}{n!}\le x^{-n}e^{H}N_n(x)\prod_{z\le p\le n}M_p(x),
\qquad N_n(x)=\sum_X\operatorname{lcm}(m,X)\prod_{\ell\in X}\bigl(e^{x^\ell/\ell}-1\bigr)
e^{H|\operatorname{pr}(X)|/z}.
\]
\textit{Here $X$ runs over the same sets as in Proposition 4.1.}

\medskip\noindent\textit{Proof.} Let $|w|=x$. We have $|e^{\varepsilon y}-1|\le e^{|y|}-1$ for
complex $y$ and $\varepsilon=\pm1$. For $|u|\le1$ we have $|G(u)|\le H$. So $|e^{G(w)}|\le e^{H}$
and $|e^{G(w^p)/p}|\le e^{H/z}$ for $p\ge z$. The function $\Gamma_p$ is entire with
$\Gamma_p(0)=1$. By the maximum principle $M_p(x)\ge1$. So for every $X$ the product of
$|\Gamma_p(w)|$ over the primes $p\notin\operatorname{pr}(X)$ in $[z,n]$ is at most
$\prod_{z\le p\le n}M_p(x)$. Hence $|\Psi_n(w)|\le e^{H}N_n(x)\prod_{z\le p\le n}M_p(x)$.
Cauchy's inequality for the coefficient of $w^n$ finishes the proof. \hfill$\square$

\medskip
In the rest of this section and in Section 5 the integers $m$ and $e$ are fixed. We take
$z=\sqrt n/\log^2n$ and $x=e^{-t}$, where $t=t_n$ satisfies $a_1\le t\sqrt{n\log n}\le a_2$ with
fixed constants $0<a_1<a_2$. Then $L=(\tfrac12+o(1))\log n$ and $1/(tL)\asymp s_n$. We take $n$ so
large that $z$ exceeds $2$ and every prime factor of $me$. We also assume $z^3>n$. A term $o(1)$ tends
to $0$ as $n\to\infty$. It may depend on $m$, $e$, $a_1$ and $a_2$.

\medskip\noindent\textbf{Lemma 4.3.} \textit{We have}
\[
\limsup_{n\to\infty}\,tL\log N_n(x)\le\nu(m,e),\qquad
\nu(m,e)=\int_0^\infty\log\Bigl(1+\sum_{j\in R}\frac{e^{-jr}}j\Bigr)dr .
\]

\noindent\textit{Proof.} Put $C=e^{H/z}$. For $\ell\le n$ let $\kappa(\ell)$ be the number of
large primes dividing $\ell$. Since $z^3>n$ we have $\kappa(\ell)\le2$. Every
$p\in\operatorname{pr}(X)$ divides some $\ell\in X$. So
$C^{|\operatorname{pr}(X)|}\le\prod_{\ell\in X}C^2$. For $y\ge0$ we have $e^y-1\le ye^y$. With
$y=x^\ell/\ell\le1/\ell$ this gives $e^{x^\ell/\ell}-1\le e^{1/\ell}x^\ell/\ell$. Put
$\Theta=eC^2$ and $\Theta_z=e^{1/z}C^2$. Put $\theta_\ell=\Theta_z$ for $\ell\ge z$ and
$\theta_\ell=\Theta$ for $\ell<z$. Also $\operatorname{lcm}(m,X)\le m\operatorname{lcm}(X)$. Hence
\begin{equation}\label{eq:Nbound}
N_n(x)\le m\sum_X\operatorname{lcm}(X)\prod_{\ell\in X}\frac{\theta_\ell x^\ell}{\ell}.
\end{equation}
Apply Lemma 2.1 to $X$. It gives a set $\mathcal D$ of pairs $(d,d')$ with $dd'\in X$ and a set
$\omega$ of remaining elements of $X$. The entries $d$ are pairwise coprime and exceed $1$. The
set $X$ is determined by $\mathcal D$ and $\omega$. By Lemma 2.1,
\[
\operatorname{lcm}(X)\prod_{\ell\in X}\frac{\theta_\ell x^\ell}{\ell}
=\prod_{(d,d')\in\mathcal D}\frac{\theta_{dd'}x^{dd'}}{d'}\prod_{\ell\in\omega}\frac{\theta_\ell x^\ell}{\ell}.
\]
So the sum in \eqref{eq:Nbound} is at most the product of a sum over $\omega$ and a sum over
$\mathcal D$. The sum over $\omega$ is at most
$\prod_{\ell\le n}(1+\Theta x^\ell/\ell)\le(1-x)^{-\Theta}$. This is $e^{O(\log n)}$, because
$1-x\ge t/2$ for small $t$. Split $\mathcal D$ into $\mathcal G_1,\mathcal G_2,\mathcal G_3$ as in
the proof of Proposition 3.2. In $\mathcal G_1$ and $\mathcal G_3$ each pair has weight at most
$\Theta/d'$. As in that proof,
$\Sigma_1\le(\Theta(n+1)(1+\log n))^z$ and $\Sigma_3\le(\Theta(n+1)(1+\log n))^{\log^4n}$. Both
are $e^{o(s_n)}$.

Let $(d,d')\in\mathcal G_2$. Then $d=p$ is a large prime. The length $pd'$ is non-free and at
least $z$. So its weight is $\Theta_zx^{pd'}/d'$. Since $pd'$ is admissible, $d'$ is coprime to
$e$. If $d'$ divided $m$ then $d'\in J$ and $pd'$ would be free. So $d'\in R$. The primes in
$\mathcal G_2$ are distinct. Hence
\[
\Sigma_2\le\prod_p\Bigl(1+\Theta_z\sum_{j\in R}\frac{x^{jp}}j\Bigr)=\exp\Bigl(\sum_pg_{\Theta_z}(pt)\Bigr),\qquad
g_a(r)=\log\Bigl(1+a\sum_{j\in R}\frac{e^{-jr}}j\Bigr).
\]
Fix $a>1$. Since $\Theta_z\to1$ we have $\Theta_z\le a$ for large $n$. The function $g_a$ is
continuous and non-increasing. We have $1\notin R$. So
$\sum_{j\in R}e^{-jr}/j\le-\log(1-e^{-r})-e^{-r}$. For $0<r\le1$ this is at most $1+\log(1/r)$.
So $g_a(r)\le\log((1+a)(2+\log(1/r)))$. For $r\ge1$ it is at most $e^{-r}$. So
$g_a(r)\le ae^{-r}$. Thus Lemma 2.2 applies with a suitable $B$. It gives
$\sum_pg_a(pt)\le(\int_0^\infty g_a+o(1))/(tL)$. Since $1/(tL)\asymp s_n$, all factors together
give $tL\log N_n(x)\le\int_0^\infty g_a+o(1)$. So the upper limit is at most
$\int_0^\infty g_a$ for every $a>1$. As $a\to1$ this integral tends to $\nu(m,e)$ by dominated
convergence. \hfill$\square$

\section{The large primes}

For $0\le r\le1$ put $\beta(r)=\max_{|u|=r}|1+G(u)|$.

\medskip\noindent\textbf{Lemma 5.1.} \textit{The function $\beta$ is continuous and non-decreasing.
We have $\beta(0)=1$ and $\beta(r)\le1+Hr$.}

\medskip\noindent\textit{Proof.} The function $1+G$ is a polynomial with value $1$ at $0$. By the
maximum principle $\beta$ is non-decreasing. Continuity is clear. For $|u|\le1$ we have
$|G(u)|\le H|u|$. \hfill$\square$

\medskip
Put $g(s)=\log\beta(e^{-s})$ for $s>0$ and
\[
I(m,e)=\int_0^\infty g(s)\,ds .
\]
By Lemma 5.1 the function $g$ is continuous, non-negative and non-increasing. It satisfies
$g(s)\le\log(1+H)$ and $g(s)\le He^{-s}$.

\medskip\noindent\textbf{Lemma 5.2.} \textit{Let $p\ge z$. Then
$\log M_p(x)\le g(pt)+H^2e^{H}e^{-2pt}/p$.}

\medskip\noindent\textit{Proof.} Let $|w|=x$ and put $u=w^p$. Then $|u|=e^{-pt}$ and
$|G(u)|\le He^{-pt}$. Write $\Gamma_p(w)=1+G(u)+E$ with $E=p(e^{G(u)/p}-1-G(u)/p)$. The series
of the exponential function gives
\[
|E|\le\frac{|G(u)|^2e^{|G(u)|/p}}{2p}\le\frac{H^2e^{H}e^{-2pt}}{p} .
\]
So $|\Gamma_p(w)|\le\beta(e^{-pt})+H^2e^{H}e^{-2pt}/p$. Finally $\log(\tau+b)\le\log\tau+b$ for
$\tau\ge1$ and $b\ge0$. \hfill$\square$

\medskip\noindent\textbf{Proposition 5.3.} \textit{We have}
\[
\sum_{z\le p\le n}\log M_p(x)\le\frac{I(m,e)+o(1)}{tL} .
\]

\noindent\textit{Proof.} The function $g$ satisfies the hypotheses of Lemma 2.2 with
$B=\max(1,H,\log(1+H)/\log2)$. So $\sum_pg(pt)=(I(m,e)+o(1))/(tL)$. Since $g\ge0$ this also bounds
the sum of $g(pt)$ over $z\le p\le n$. By Lemma 5.2 it remains to bound $\sum_pe^{-2pt}/p$. This
sum is at most $\sum_{k\le1/t}1/k+t\sum_{k>1/t}e^{-2kt}=O(L)$. This is $o(1/(tL))$.
\hfill$\square$

\medskip\noindent\textbf{Proposition 5.4.} \textit{For fixed $m$ and $e$,}
\[
\log\frac{|S_n(m,e)|}{n!}\le\Bigl(2\sqrt{2\bigl(I(m,e)+\nu(m,e)\bigr)}+o(1)\Bigr)s_n .
\]

\noindent\textit{Proof.} Put $A'=I(m,e)+\nu(m,e)$. We have $A'>0$, because $1+G$ is not constant.
Take $t=\sqrt{2A'/(n\log n)}$. It satisfies the condition of Section 4 for any fixed
$a_1<\sqrt{2A'}<a_2$. Then $nt=\sqrt{2A'}\,s_n$ and
$1/(tL)=(2+o(1))s_n/\sqrt{2A'}$. Lemma 4.2, Lemma 4.3 and Proposition 5.3 give
\[
\log\frac{|S_n(m,e)|}{n!}\le nt+H+\frac{A'+o(1)}{tL}=(2\sqrt{2A'}+o(1))s_n .
\]
\hfill$\square$

\section{Proof of Theorem 1.1(b) and of the corollaries}

\noindent\textbf{Lemma 6.1.} \textit{Let $Q$ be a finite set of odd primes. For $T\subseteq Q$ put
$m_T=2\prod_{p\in Q\setminus T}p$ and $e_T=\prod_{p\in T}p$. Then}
\[
D_n^*=\sum_{T\subseteq Q}\Bigl(\prod_{p\in T}(1-p)\Bigr)S_n(m_T,e_T).
\]

\noindent\textit{Proof.} Let $\sigma\in S_n$ and $o=o(\sigma)$. For $p\in Q$ put $\delta_p=1$ if
$p\nmid o$ and $\delta_p=0$ otherwise. The primes in $Q$ are odd and distinct. So
$\operatorname{lcm}(m_T,o)=\operatorname{lcm}(2,o)\prod_{p\in Q\setminus T}p^{\delta_p}$. Also
$\gcd(o,e_T)=1$ exactly when $\delta_p=1$ for every $p\in T$. Hence $\sigma$ contributes
\[
\operatorname{sgn}(\sigma)\operatorname{lcm}(2,o)\prod_{p\in Q}\bigl(p^{\delta_p}+(1-p)\delta_p\bigr)
\]
to the right side. Each factor of the product equals $1$. Summing over $\sigma$ gives $D_n^*$.
\hfill$\square$

\medskip
For $Q=\emptyset$ the lemma is trivial. Proposition 5.4 then applies directly to $D_n^*=S_n(2,1)$.
Numerically $I(2,1)+\nu(2,1)=0.9772\ldots$ This lies between $A_2$ and $A$. So the direct route gives
Proposition 6.3 below. For Theorem 1.1(b) we take a larger set $Q$. It makes more multipliers free. This lowers $\nu$. The integral $I$ grows slightly. For
$Q=\{3\}$ the lemma reads
\[
D_n^*=S_n(6,1)-2S_n(2,3).
\]
For both pairs the value of $I+\nu$ is below $A_2$ (Proposition 6.2(c)). So $Q=\{3\}$ already proves
Theorem 1.1(b), with $\delta=0.04$. We use $Q=\{3,5\}$ because it gives a better rate. Then
\[
D_n^*=S_n(30,1)-2S_n(10,3)-4S_n(6,5)+8S_n(2,15).
\]
For these four pairs the sets $J$ are $\{1,2,3,5,6,10,15,30\}$, $\{1,2,5,10\}$, $\{1,2,3,6\}$ and
$\{1,2\}$. For larger sets a floating-point computation gives the largest value of $I+\nu$ as
$0.8729\ldots$ for $Q=\{3,5,7\}$ and as $0.8692\ldots$ for $Q=\{3,5,7,11\}$. These two values serve
only for comparison. The proofs do not use them.

\medskip\noindent\textbf{Proposition 6.2.} \textit{{\rm(a)} We have $0.9356\le A_2\le0.9357$ and
$1.1178\le A\le1.1179$.}

\textit{{\rm(b)} We have}
\[
I(30,1)+\nu(30,1)\le0.8839,\qquad I(10,3)+\nu(10,3)\le0.8451,
\]
\[
I(6,5)+\nu(6,5)\le0.8653,\qquad I(2,15)+\nu(2,15)\le0.8263 .
\]

\textit{{\rm(c)} We have}
\[
I(6,1)+\nu(6,1)\le0.9041,\qquad I(2,3)+\nu(2,3)\le0.8642,\qquad I(2,1)+\nu(2,1)\le0.9777 .
\]

\noindent\textit{Proof.} The proof is computer-assisted. All integrands are non-negative and
non-increasing. For $I(m,e)$ this is part of Lemma 5.1. So a Riemann sum with left end points is an
upper bound for the integral over the range of the nodes. A Riemann sum with right end points is a
lower bound. All nodes are binary floating point numbers. They enter the computation exactly. Every
function value and every sum is computed in ball arithmetic with at least $160$ bits, using the
library Arb \cite{Johansson2017}. The reported bounds are outer end points of the resulting balls.

For $A$, $A_2$ and the integrals $\nu(m,e)$ the nodes form a geometric progression with ratio
$1+2^{-16}$ from $2^{-40}$ to $1$. Beyond $1$ the step is $2^{-16}$ up to $40$. These are $4372963$
nodes. For $\nu(m,e)$ we use
\[
\sum_{j\in R}\frac{x^j}{j}=-\sum_{d\mid e}\frac{\mu(d)}{d}\log(1-x^d)-\sum_{j\in J}\frac{x^j}{j},
\]
where $\mu$ is the M\"obius function. Each integrand is at most $-\log(1-e^{-r})\le r+\log(1/r)$ for
$0<r<2^{-40}$. For $r>40$ the integrands of $A$ and $A_2$ are at most $e^{-r}/(1-e^{-40})$. Let $j_0$
be the least element of $R$. For $r>40$ the integrand of $\nu(m,e)$ is at most
$e^{-j_0r}/(j_0(1-e^{-40}))$. We have $j_0=3$ for the pair $(2,1)$ and $j_0=4$ for the other six
pairs. The integrals of these bounds over $(0,2^{-40})$ and over $(40,\infty)$ are added to the upper
bounds. Together they are less than $10^{-10}$. For the lower bounds the two ranges are dropped.

For $I(m,e)$ the nodes have step $2^{-10}$ on $[0,3]$ and step $2^{-8}$ on $[3,12]$. These are
$3072+2304=5376$ nodes. For $s>12$ we use $g(s)\le He^{-s}$. So the term $He^{-12}$ is added. At a
node $s$ the number $r=e^{-s}$ is enclosed in a ball. Put $P(\alpha)=|1+G(re^{i\alpha})|^2$. Then
$P(\alpha)=C_0+\sum_{d=1}^mC_d\cos(d\alpha)$ with explicit polynomials $C_d$ in $r$. The function
$P$ is even. We have $|P''|\le K$ with $K=\sum_dd^2|C_d|$. The numbers $C_d$ and $K$ are computed as
balls from the ball for $r$. Let $\Delta=\pi/4096$. On a cell
$[i\Delta,(i+1)\Delta]$ the function $P$ is at most the larger of its two end point values plus
$\Delta^2K/8$. Hence
\[
\beta(r)^2\le\max_{0\le i\le4096}P(i\Delta)+\frac{\Delta^2K}{8} .
\]
The values $\cos(di\Delta)$ are enclosed in balls. The upper end point of the right side is a bound
for $\beta(r)^2$. The computation gives
\[
0.935643\le A_2\le0.935653,\qquad1.117858\le A\le1.117870,
\]
\[
I(30,1)\le0.739171,\quad I(10,3)\le0.719636,\quad I(6,5)\le0.724472,\quad I(2,15)\le0.704624,
\]
\[
\nu(30,1)\le0.144646,\quad \nu(10,3)\le0.125425,\quad \nu(6,5)\le0.140820,\quad
\nu(2,15)\le0.121589 .
\]
The pairs $(2,3)$ and $(2,1)$ have the same set $J=\{1,2\}$ as $(2,15)$. So
$I(2,3)=I(2,1)=I(2,15)$. The pair $(6,1)$ has the same set $J$ as $(6,5)$. So $I(6,1)=I(6,5)$. The
computation gives
\[
\nu(6,1)\le0.179559,\qquad \nu(2,3)\le0.159549,\qquad \nu(2,1)\le0.273005 .
\]
The end points of the balls for $A$ and $A_2$ give the sharper bounds
\begin{equation}\label{eq:cbounds}
2.735900\le c_2\le2.735914,\qquad0.25454<c-c_2<0.25458 .
\end{equation}
A second computation uses the interval arithmetic of \cite{mpmath} with coarser nodes. It gives
$A_2\ge0.9353$. For the four pairs of part (b) it gives $I+\nu\le0.8863$. Here
$2\sqrt{2\cdot0.9353}-2\sqrt{2\cdot0.8863}>0.07$. So the second computation confirms the margin used in
the proof of Theorem 1.1(b). For the three pairs of part (c) it gives the bounds $0.9062$,
$0.8663$ and $0.9798$. They suffice for the two uses of part (c).
\hfill$\square$

\medskip\noindent\textit{Proof of Theorem 1.1(b).} By Lemma 3.1 and Lemma 6.1 with $Q=\{3,5\}$,
\[
D_n+V_n=S_n(30,1)-2S_n(10,3)-4S_n(6,5)+8S_n(2,15).
\]
By Proposition 5.4 and Proposition 6.2(b) the four sums on the right are at most
$n!\exp((c^*+o(1))s_n)$ in absolute value, where $c^*=2\sqrt{2\cdot0.8839}<2.6592$. By Theorem
1.1(a) and Proposition 6.2(a) we have $V_n\ge n!\exp((c_2-o(1))s_n)$ with
$c_2\ge2\sqrt{2\cdot0.9356}>2.7358$. Hence
\[
|D_n+V_n|\le15V_n\exp\bigl(-(0.0766-o(1))s_n\bigr)\le V_n\exp(-0.07\,s_n)
\]
for all sufficiently large $n$. This proves the theorem with $\delta=0.07$. With $Q=\{3\}$ the same
argument uses Proposition 6.2(c). There $2\sqrt{2\cdot0.9041}<2.6894$. This gives $\delta=0.04$.
\hfill$\square$

\medskip\noindent\textbf{Proposition 6.3.} \textit{For all sufficiently large $n$ we have
$|D_n|\le T_n\,e^{-0.19\,s_n}$. Hence $D_n=o(T_n)$ and $\psi(S_n)/\psi(A_n)\to2$. In particular
$3\psi(A_n)>\psi(S_n)$ for all sufficiently large $n$.}

\medskip\noindent\textit{Proof.} By Lemma 3.1 we have $|D_n|\le|D_n^*|+V_n$. Here $D_n^*=S_n(2,1)$.
By Proposition 5.4 and Proposition 6.2(c) we have $|D_n^*|\le n!\exp((2.7968+o(1))s_n)$, since
$2\sqrt{2\cdot0.9777}<2.7968$. By Proposition 3.2 and \eqref{eq:cbounds} we have
$V_n\le n!\exp((2.7360+o(1))s_n)$. By the theorem of Goh and Schmutz \cite{GohSchmutz1991} and
Proposition 6.2(a) we have $T_n=n!\exp((c+o(1))s_n)$ with $c\ge2\sqrt{2\cdot1.1178}>2.9903$. Since
$2.9903-2.7968>0.19$, the first claim follows. The limit follows from \eqref{eq:ratio}. The
inequality $3\psi(A_n)>\psi(S_n)$ says that $D_n>-T_n/3$. \hfill$\square$

\medskip\noindent\textit{Proof of Corollary 1.2.} Part (a) follows from Theorem 1.1(b). By Theorem
1.1 we have $\log(|D_n|/n!)=\log(V_n/n!)+o(1)=(c_2+o(1))s_n$. The theorem of Goh and Schmutz
\cite{GohSchmutz1991} states that $\log(T_n/n!)=(c+o(1))s_n$. Together this gives
$\log(|D_n|/T_n)=-(c-c_2+o(1))s_n$. The bounds for $c-c_2$ are in \eqref{eq:cbounds}.
\hfill$\square$

\medskip\noindent\textit{Proof of Corollary 1.3.} By \eqref{eq:ratio} and $D_n/T_n\to0$ we have
$\psi(S_n)/\psi(A_n)=2-(2+o(1))D_n/T_n$. By Theorem 1.1(b) we have $D_n=-(1+o(1))V_n$. Since
$V_n>0$, the ratio exceeds $2$ for all sufficiently large $n$. \hfill$\square$

\section{Exact computations and a conjecture}

The order of a permutation depends only on the set of its cycle lengths. We computed $T_n$
and $D_n$ for $n\le150$ by a dynamic programme over the pairs (number of points used, least
common multiple so far). Cycle lengths $2,3,\dots,n$ are added in turn with all multiplicities.
Signs are carried along. Fixed points are attached at the end by a binomial convolution. The
same programme restricted to odd cycle lengths gives $V_n$. All arithmetic is in exact integers.
For $n\le50$ a direct sum over cycle types confirms these values and Lemma 3.1. It also confirms
Lemma 6.1 for $Q=\{3\}$, $\{3,5\}$ and $\{3,5,7\}$ in the ranges $n\le50$, $n\le55$ and $n\le30$
respectively. Proposition 4.1 was checked in exact rational arithmetic for $n\le15$ and the seven
pairs $(m,e)$ of Proposition 6.2.

\begin{table}[ht]
\centering
\begin{tabular}{@{}rrrr@{}}
\toprule
$n$ & $\psi(A_n)$ & $\psi(S_n)$ & $V_n$\\
\midrule
3 & 7 & 13 & 7\\
4 & 31 & 67 & 25\\
5 & 211 & 471 & 181\\
6 & 1411 & 3271 & 961\\
7 & 12601 & 31333 & 8611\\
8 & 137047 & 299223 & 91057\\
9 & 1516831 & 3291487 & 939625\\
10 & 18111751 & 39020911 & 10163041\\
11 & 223179001 & 543960561 & 119233951\\
12 & 2973194071 & 7466726983 & 1567163401\\
13 & 46287964867 & 118551513523 & 21299079517\\
\bottomrule
\end{tabular}
\caption{Exact values for $3\le n\le13$.}
\label{tab:small}
\end{table}

\begin{table}[ht]
\centering
\begin{tabular}{@{}rrrr@{}}
\toprule
$n$ & $\psi(S_n)/\psi(A_n)$ & $V_n/T_n$ & $-D_n/V_n$\\
\midrule
10 & 2.15445 & 0.26045 & 0.275\\
20 & 2.38300 & 0.15001 & 1.071\\
30 & 2.43692 & 0.10287 & 1.743\\
40 & 2.15814 & 0.10180 & 0.720\\
50 & 2.18244 & 0.07698 & 1.086\\
60 & 2.24686 & 0.06331 & 1.735\\
70 & 2.11465 & 0.06131 & 0.884\\
80 & 2.13610 & 0.05208 & 1.223\\
90 & 2.15585 & 0.04518 & 1.600\\
100 & 2.08335 & 0.04439 & 0.901\\
110 & 2.09698 & 0.03880 & 1.192\\
120 & 2.10990 & 0.03502 & 1.488\\
130 & 2.06686 & 0.03397 & 0.952\\
140 & 2.07679 & 0.03053 & 1.211\\
150 & 2.08371 & 0.02810 & 1.429\\
\bottomrule
\end{tabular}
\caption{Ratios for selected $n$.}
\label{tab:ratios}
\end{table}

\medskip\noindent\textit{Proof of Proposition 1.4.} The exact values give $D_n>-T_n/3$ for
$3\le n\le150$. Among these $n$ the ratio $\psi(S_n)/\psi(A_n)$ is largest at $n=13$. The values
give $D_n<0$ for $n=2$ and $4\le n\le150$. \hfill$\square$

\medskip
Table~\ref{tab:ratios} shows that $-D_n/V_n$ oscillates around $1$. For $100\le n\le150$ it lies
between $0.7255$ and $1.4942$. The values for $n\le150$ belong to a pre-asymptotic range. At
$n=150$ we have $s_n=5.47$, $\log(V_n/n!)/s_n=1.68$ and $\log(T_n/n!)/s_n=2.33$. The limits are
$c_2=2.7359\ldots$ and $c=2.9904\ldots$ If the terms $o(1)$ are ignored, the bound $15e^{-0.0766s_n}$
from the proof of Theorem 1.1(b) falls below $0.01$ for $s_n>95.47$, that is for $n>1.054\cdot10^5$.

Explicit values of the thresholds in our results need three further inputs. The first is the prime
number theorem with an explicit error term in Lemma 2.2, for instance through the bounds of Rosser and
Schoenfeld \cite{RosserSchoenfeld}. This lemma enters Propositions 3.2, 3.3 and 5.3. It also enters
Lemma 4.3. The second is a rate of convergence in Lemma 2.4 together with quantitative forms of the
limits $y\to\infty$ in Proposition 3.3 and $a\to1$ in Lemma 4.3. The third is an explicit constant in
the error term of \eqref{eq:GS} \cite{Stong1998}. This third input concerns only Corollary 1.2(b),
Corollary 1.3 and Proposition 6.3.

\medskip\noindent\textbf{Conjecture 7.1.} \textit{For every $n\ge3$ we have $3\psi(A_n)>\psi(S_n)$.
For every $n\ge4$ we have $D_n<0$.}

\medskip
By Proposition 1.4, Proposition 6.3 and Corollary 1.2 both statements hold for $n\le150$ and for all
sufficiently large $n$. The same is true of the sharper inequality
$\psi(S_n)/\psi(A_n)\le\psi(S_{13})/\psi(A_{13})$ for $n\ge3$ considered in \cite{Singh2026}.

\section*{Data availability}

The programs and data that support this study are available from the author on request. They
reproduce Proposition 6.2, all the checks listed in Section 7 and the two tables.

\section*{Competing interests}

The author has no competing interests to declare. No funding was received for this work.

\end{document}